\documentclass[12pt,twoside]{article}
\usepackage[hmargin=0.8in,vmargin=0.9in]{geometry}
\usepackage{fancyhdr}
\usepackage{graphicx}
\usepackage{subfigure}
\usepackage{amssymb}
\usepackage{amsmath}
\usepackage{amsfonts}
\usepackage{theorem}
\usepackage{mathrsfs}
\usepackage{mathtools}
\usepackage{bm}
\usepackage{color}
\usepackage{setspace}
\usepackage{exscale}
\usepackage{relsize}
\usepackage{epstopdf}
\usepackage{float}
\usepackage{cite}
\usepackage{cleveref}
\usepackage{refcheck}
\usepackage{nicefrac}
\usepackage{setspace}
\usepackage{booktabs,multirow} % for much better looking tables
\usepackage{array} % for better arrays (eg matrices) in maths
\usepackage{paralist} % very flexible & customisable lists (eg. enumerate/itemize, etc.)
\usepackage{verbatim} % adds environment for commenting out blocks of text & for better verbatim
\usepackage{subfigure} % make it possible to include more than one captioned figure/table in a single float
\theoremstyle{plain}			% use "default" font
\newtheorem{thm}{Theorem}[section]

\newtheorem{conj}[thm]{Conjecture}

\newtheorem{rmk}[thm]{Remark}
{\theorembodyfont{\rmfamily}}
\newenvironment{acknowledgment}{{\flushleft \bf Acknowledgment:}}{}

\usepackage{tikz}
\usetikzlibrary{positioning}
\usepackage{xcolor}

\allowdisplaybreaks[1]

\numberwithin{equation}{section}
\numberwithin{figure}{section}
\numberwithin{table}{section}

\newcommand\eref[1]{(\ref{#1})}

\newcommand*\xbar[1]{%
  \hbox{%
    \vbox{%
      \hrule height 0.5pt % The actual bar
      \kern0.4ex%         % Distance between bar and symbol
      \hbox{%
        \kern-0.05em%      % Shortening on the left side
        \ensuremath{#1}%
        \kern-0.00em%      % Shortening on the right side
      }%
    }%
  }%
}

\newcommand{\mo}{\bm{0}}

\newcommand{\dx}{\Delta x}
\newcommand{\dy}{\Delta y}

\def\softd{{\leavevmode\setbox1=\hbox{d}%
          \hbox to 1.05\wd1{d\kern-0.4ex{\char039}\hss}}}

\title{On the Gelfand Problem and Viscosity Matrices for Two-Dimensional Hyperbolic Systems of Conservation Laws}
\author{Shaoshuai Chu\thanks{Department of Mathematics and Shenzhen International Center for Mathematics, Southern University of Science and
Technology, Shenzhen, 518055, China; {\tt chuss2019@mail.sustech.edu.cn}}, ~Igor Kliakhandler\thanks{Department of Mathematics, Michigan
Technological University, Houghton, MI 49931, USA; {\tt igor@mtu.edu}}, ~and Alexander Kurganov\thanks{Department of Mathematics, Shenzhen
International Center for Mathematics, and Guangdong Provincial Key Laboratory of Computational Science and Material Design, Southern
University of Science and Technology, Shenzhen, 518055, China; {\tt alexander@sustech.edu.cn}}}
\begin{document}

\date{}
\maketitle
\begin{abstract}
We present counter-intuitive examples of a viscous regularizations of a two-dimensional strictly hyperbolic system of conservation laws. The
regularizations are obtained using two different viscosity matrices. While for both of the constructed ``viscous'' systems waves propagating
in either $x$- or $y$-directions are stable, oblique waves may be linearly unstable. Numerical simulations fully corroborate these
analytical results. To the best of our knowledge, this is the first nontrivial result related to the multidimensional Gelfand problem. 
Our conjectures provide direct answer to Gelfand's problem both in one- and multi-dimensional cases. 
\end{abstract}



\section{Introduction}
Hyperbolic systems of conservation laws and their viscous regularizations pertaining to the natural sciences are among the most important
partial differential equations (PDEs) and invariably stir the interest of mathematicians of all walks. Stability of solutions and their
dependence on parameters are among the basic and most important questions that need to be investigated. During the height of the Cold War,
Soviet-American mathematician I. M. Gelfand published a paper with an unassuming title ``Some problems in the theory of quasi-linear
equations'' \cite{Gelfand59}. It was the time when high-temperature physics, nuclear reactions, and powerful explosions became the hotly
pursued and studied areas as tools of political wrestling. Being inspired by the mathematical problems relevant to the new realms of
physics, with his (un)usual insight and tour de force, Gelfand formulated and posted a few fundamental problems in the theory of PDEs.

Among them, one of the most intriguing questions relates to the role and structure of dissipative, ``viscous'' effects, which are usually
present as higher (usually second) order partial derivative terms in viscous regularizations of the hyperbolic conservation laws. Being
typically small, those dissipative effects become important near shock waves or areas with sharp spatial gradients. The question asked by
Gelfand, sounds deceptively simple: What structure of dissipative terms will ensure the convergence of solutions of the regularized
hyperbolic conservation systems to those of the corresponding inviscid ones?

We consider strictly hyperbolic systems, whose Jacobians have distinct real eigenvalues. For those systems an (obvious) requirement for the
viscosity matrix is to be positive definitive, that is, ``viscous''. Of course, the diagonal viscosity matrices are the simplest, and 
correspond to the viscous terms in the Navier-Stokes equations of fluid motion. A diagonal viscous matrix with positive entries was used by
Bianchini and Bressan \cite{BB05} to prove the global existence and uniqueness for one-dimensional (1-D) strictly hyperbolic systems via the
vanishing viscosity method. Nevertheless, Majda and Pego \cite{MP85} discovered a remarkable phenomenon: viscous shocks for 1-D systems of
convection-diffusion systems with strictly hyperbolic convective term and small and positive definitive, but {\em non}-diagonal viscosity
matrix may {\em not} converge to the shock waves of the corresponding inviscid system. Moreover, such system may be linearly unstable. 
Those instability examples have direct pertinence to the modeling of natural phenomena in larger, ``macro'' scales when both viscous and
convective terms are not expected to be simply diagonal. Later, Kliakhandler \cite{Kliakhandler99} found the same phenomenon for specific
fluid mechanics problems.

In \cite{MP85}, Majda and Pego considered the simplest 1-D linear system with constant coefficients,
\begin{equation}
\bm u_t+A\bm u_x=C\bm u_{xx},
\label{1.1}
\end{equation}
where $x$ is a spatial variable, $t$ is time, $\bm u\in\mathbb R^N$ is the vector of unknowns ($N\ge2$), the matrix $A$ has distinct real
eigenvalues, and the matrix $C$ has real positive eigenvalues. Even though one may expect the solutions of \eref{1.1} to be smooth and
stable, Majda and Pego observed that the simplest solution, $\bm u(x,t)=e^{ikx+\omega t}\hat{\bm u}$ with $\hat{\bm u}={\bf Const}$ and the
linear stability operator $\Omega(k):=-ikA-k^2C$ may have nontrivial behavior. It is obvious that for small $k$ (the so-called
long-wavelength solutions), the convective part is dominating and, in particular, the wave speeds are defined by the eigenvalues of $A$.
This is similar to the case of a scalar convection-diffusion equation. However, unlike the scalar case, the stability of the system
\eref{1.1} is ensured by the negativity of the real parts of the eigenvalues of the matrix $\Omega$, which is not guaranteed by the
positivity of the eigenvalues of $C$, which dominates the eigenvalue problem for large $k$ only. The behavior of the eigenvalues of
$\Omega(k)$ for small $k$ may be understood by the perturbative analysis of the simple ``wavy'' eigenvalues of the matrix $A$. Those
perturbative terms include interplay between the matrices $A$ and $C$. It turns out that $\Omega$ may have eigenvalues with positive real
part for small wavenumbers $k$ so that the corresponding solutions of \eref{1.1} will grow indefinitely in time, contrary to the ``naive''
expectations. This may happen when $C$ is positive definitive but not diagonal. We stress that the non-diagonality of both the advective
matrix $A$ and the viscosity matrix $C$ is critical for the linear instability.

Concerning the multidimensional systems, very little has been known about the influence of viscous terms on stability and convergence. In
this paper, we discover the two-dimensional (2-D) phenomenon that extends, in a certain sense, the 1-D results of Majda and Pego
\cite{MP85}. To this end, we consider the 2-D extension of \eref{1.1},
\begin{equation}
\bm u_t+A\bm u_x+B\bm u_y=C\bm u_{xx}+D\bm u_{yy},
\label{1.2}
\end{equation}
where the advective matrices $A$ and $B$ have distinct real eigenvalues and the viscosity matrices $C$ and $D$ are positive definite, and
show that viscous regularizations of 2-D strictly hyperbolic systems may be stable in the directions of both $x$- and $y$-axes, but yet some
oblique waves may be unstable. We present an explicit example of such system (linearized shallow water equations), for which we perform
linear stability analysis. We also conduct numerical simulations for this unstable linear system and also for nonlinear regularized
Saint-Venant system. Our numerical results fully corroborate the obtained linear stability results. Finally, we conjecture on the stability
of the eigenvalue problem and point to the obvious three-dimensional (3-D) extensions of our results.

\section{Two-Dimensional Saint-Venant Systems}
The Saint-Venant system introduced in \cite{Saint71} is one of the most commonly used models of shallow water flows in rivers or coastal
areas.

In the case of flat bottom topography, the 2-D Saint-Venant system reads as
\begin{equation}
\begin{cases}
h_t+q_x+p_y=0,\\[0.5ex]
q_t+\left(\dfrac{q^2}{h}+\dfrac{g}{2}h^2\right)_x+\left(\dfrac{pq}{h}\right)_y=0,\\[2.5ex]
p_t+\left(\dfrac{pq}{h}\right)_x+\left(\dfrac{p^2}{h}+\dfrac{g}{2}h^2\right)_y=0.
\end{cases}
\label{2.1}
\end{equation}
Here, $h$ is the water depth, $u$ and $v$ are the $x$- and $y$-velocities, respectively, $q=hu$ and $p=hv$ are the corresponding discharges,
and $g$ is the acceleration due to gravity.

The 2-D Saint-Venant system \eref{2.1} can be linearized about a ``lake-at-rest'' steady state $h(x,y,t)\equiv h_0=\texttt{Const}$,
$u(x,y,t)=v(x,y,t)\equiv0$. To this end, we introduce the perturbation variables $h'$, $q'$, and $p'$ such that $h=h_0+h'$, $q=q'$, and
$p=p'$, and substitute them into  \eref{2.1} to obtain
\begin{equation}
\begin{cases}
h'_t+q'_x+p'_y=0,\\[0.5ex]
q'_t+\left(\dfrac{(q')^2}{h_0+h'}+\dfrac{g}{2}(h_0+h')^2\right)_x+\left(\dfrac{p'q'}{h_0+h'}\right)_y=0,\\[2.5ex]
p'_t+\left(\dfrac{p'q'}{h_0+h'}\right)_x+\left(\dfrac{(p')^2}{h_0+h'}+\dfrac{g}{2}(h_0+h')^2\right)_y=0.
\end{cases}
\label{2.2}
\end{equation}
We then neglect the nonlinear terms in \eref{2.2} and end up with the 2-D linearized Saint-Venant system:
\begin{equation*}
\begin{cases}
h'_t+q'_x+p'_y=0,\\
q'_t+gh_0h'_x=0,\\
p'_t+gh_0h'_y=0,
\end{cases}
\end{equation*}
which can be rewritten in the vector form
\begin{equation}
\bm u_t+A\bm u_x+B\bm u_y=\mo,
\label{2.3}
\end{equation}
with
\begin{equation}
\bm u=\begin{pmatrix}h'\\q'\\p'\end{pmatrix},\quad A=\begin{pmatrix}0&1&0\\gh_0&0&0\\0&0&0\end{pmatrix},\quad
B=\begin{pmatrix}0&0&1\\0&0&0\\gh_0&0&0\end{pmatrix}.
\label{2.4}
\end{equation}

\section{Viscous Approximations of the Saint-Venant Systems}\label{sec3}
We now add the viscosity terms to the right-hand sides (RHS) of \eref{2.3}--\eref{2.4} and study the stability of the resulting viscous
linearized shallow water equations, which read as \eref{1.2}, \eref{2.4}. We consider the following specific viscosity matrices $C$ and $D$:
\begin{equation}
C=\begin{pmatrix}\varepsilon&0&0\\0&\varepsilon&0\\0&0&\varepsilon\end{pmatrix},\quad
D=\begin{pmatrix}\varepsilon&0&0\\0&\varepsilon&0\\\Delta&\Delta&\varepsilon\end{pmatrix},
\label{3.1}
\end{equation}
where $\varepsilon>0$ and $\Delta$ are constants.

We now analyze the stability of the linearized viscous shallow water equations \eref{1.2}, \eref{2.4}, \eref{3.1}. To this end, we take the
following ansatz:
\begin{equation}
\begin{pmatrix}h'\\q'\\p'\end{pmatrix}=\begin{pmatrix}\hat h_1'\\\hat q_1^{\,\prime}\\\hat p_1^{\,\prime}\end{pmatrix}
e^{\omega t+ik\left(\gamma x+\sqrt{1-\gamma^2}y\right)},
\label{3.2}
\end{equation}
where $\omega$ is a wave magnitude, $k$ is its frequency, $\gamma$ stands for a spatial direction of the wave, and
$(\hat h_1',\hat q_1^{\,\prime},\hat p_1^{\,\prime})^\top$ is a constant nonzero vector. Substituting \eref{3.2} into \eref{1.2},
\eref{2.4}, \eref{3.1} results in
\begin{equation*}
E\bm u=\mo,
\end{equation*}
where
\begin{equation*}
E=\begin{pmatrix}
\omega+\varepsilon k^2\gamma^2+\varepsilon k^2(1-\gamma^2)&ik\gamma&ik\sqrt{1-\gamma^2}\\
ik\gamma gh_0&\omega+\varepsilon k^2\gamma^2+\varepsilon k^2(1-\gamma^2)&0\\
ik\sqrt{1-\gamma^2}gh_0+\Delta k^2(1-\gamma^2)&\Delta k^2(1-\gamma^2)&\omega+\varepsilon k^2\gamma^2+\varepsilon k^2(1-\gamma^2)
\end{pmatrix}.
\end{equation*}
Therefore, $\bm u$ is an eigenvector that corresponds to the zero eigenvalue if and only if the following characteristic equation is
satisfied:
\begin{equation*}
\begin{aligned}
\varepsilon^3k^6-\Delta k^3(1-\gamma^2)^{\frac{3}{2}}(gh_0k\gamma+i\omega)+\varepsilon k^2(gh_0k^2-i\Delta k^3(1-\gamma^2)^{\frac{3}{2}}
&+3\omega^2)+3\varepsilon^2k^4\omega\\
&+(gh_0k^2+\omega^2)\omega=0.
\end{aligned}
\end{equation*}
Its solution $\omega$ depends on the wave number $k$ and the parameters $\varepsilon$, $g$, $h_0$, and $\Delta$. Clearly, if $\omega$ is
positive for some $k$, then the system \eref{1.2}, \eref{2.4}, \eref{3.1} is unstable.

\begin{figure}[ht!]
\centerline{\includegraphics[trim=0.cm 0.cm 0.0cm 0.cm, clip, width=7.cm]{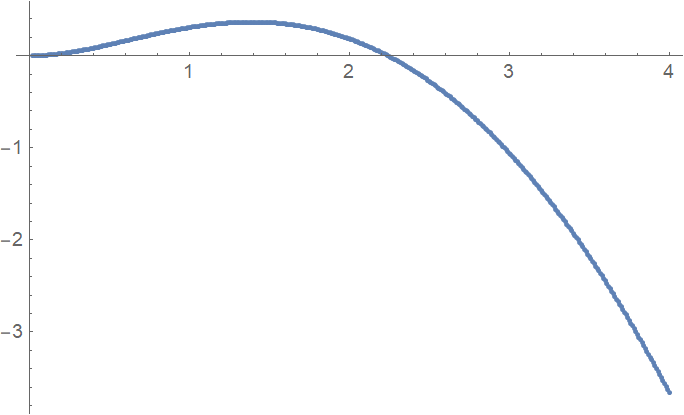} \hspace{1cm}
            \includegraphics[trim=0.cm 0.cm 0.0cm 0.cm, clip, width=7.cm]{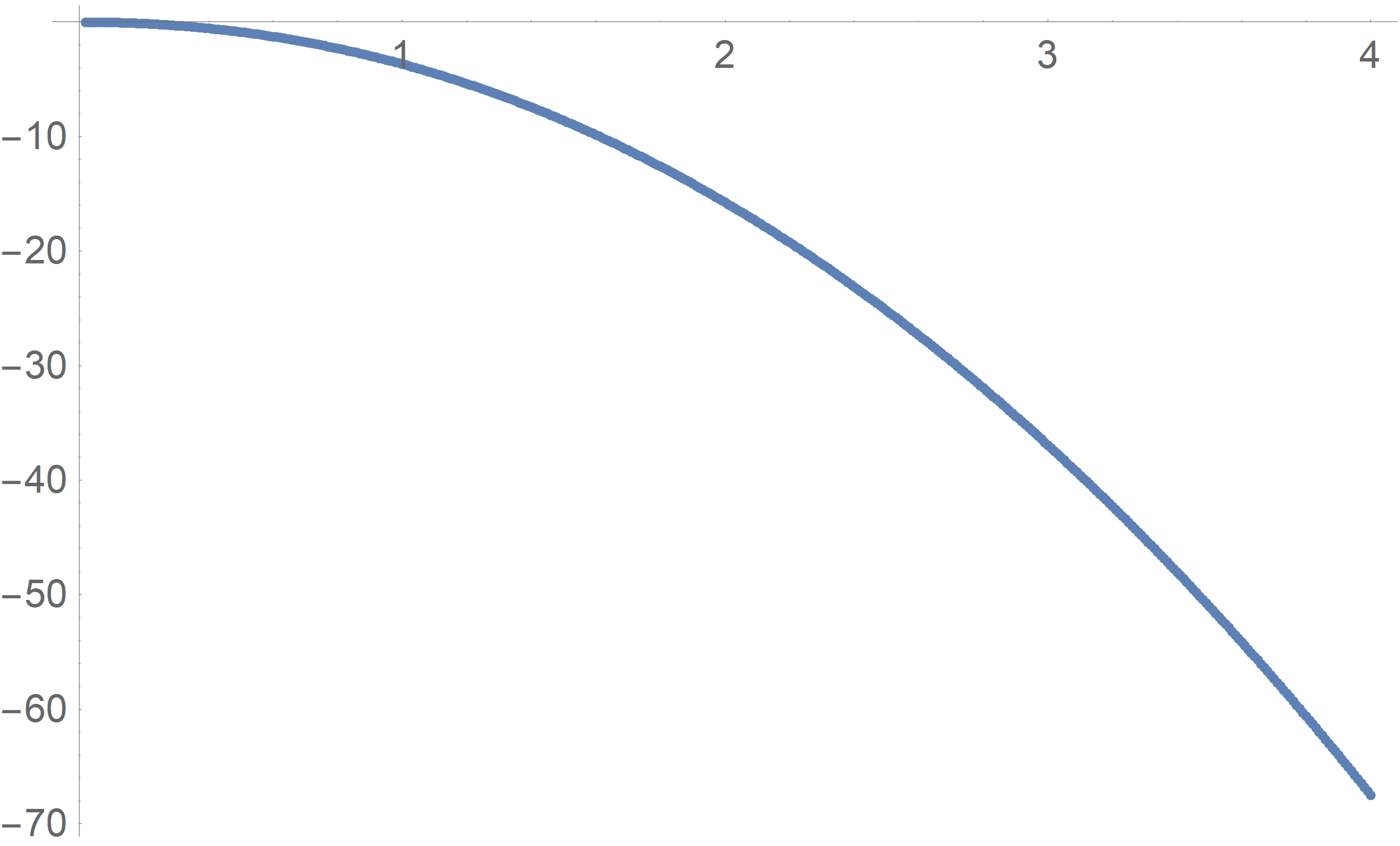}}
\caption{\sf $w$ as a function of $k$ for $\varepsilon=1$ (left) and $\varepsilon=5$ (right).\label{fig1}}
\end{figure}

Next, we take following parameter values: $h_0=1$, $g=10$, $\Delta=5$, and $\gamma=1/2$, and then plot the function $\omega(k)$ for
$\varepsilon=1$ and $\varepsilon=5$ in Figure \ref{fig1}. As one can see, when $\varepsilon=5$, $\omega(k)\le0$ for all $k$ and thus the
system is stable. However, when $\varepsilon=1$, $\omega(k)>0$ for $0<k\lesssim2.3$ with the maximum value about 0.4 achieved at about
$k=1.5$. This means that $\varepsilon=1$ corresponds to the unstable regime.

\section{Numerical Examples}
In this section, we conduct several numerical experiments using a semi-discrete fifth-order finite-difference alternative weighted
essentially non-oscillatory (A-WENO) scheme described in \cite[\S3]{CKN21} (see also \cite[Appendix C]{CKKO23}). The diffusion terms have
been discretized using a standard six-order finite-difference approximations of $\bm u_{xx}$ and $\bm u_{yy}$. We have used the three-stage
third-order strong stability preserving (SSP) Runge-Kutta method (see, e.g., \cite{Gottlieb12,Gottlieb11}) for the temporal discetization.

In all of the numerical examples, periodic boundary conditions are imposed on the four sides of the computational domain
$[-30,30]\times[-30,30]$, which corresponds to 15 wavelengths in both the $x$- and $y$-directions. We compute the numerical solutions until
a very large final time $t=250$ on a uniform mesh with $\dx=\dy=0.2$.

\subsection{Linear Case}
We first conduct numerical experiments for the viscous linearized shallow water equations \eref{1.2}, \eref{2.4}, \eref{3.1} subject to the
following initial conditions:
\begin{equation*}
h'(x,y,0)=\begin{cases}10^{-2}\Big(1+\dfrac{1}{40}e^{2(1-x^2-y^2)}\Big)&\mbox{if }x^2+y^2<1,\\10^{-2},&\mbox{otherwise},\end{cases}\quad
q'(x,y,0)=p'(x,y,0)\equiv0,
\end{equation*}
which are a perturbation of the steady state $(h',q',p')\equiv(10^{-2},0,0)$.

\begin{figure}[ht!]
\centerline{\includegraphics[trim=0.0cm 0.3cm 0.3cm 0.2cm, clip, width=7.cm]{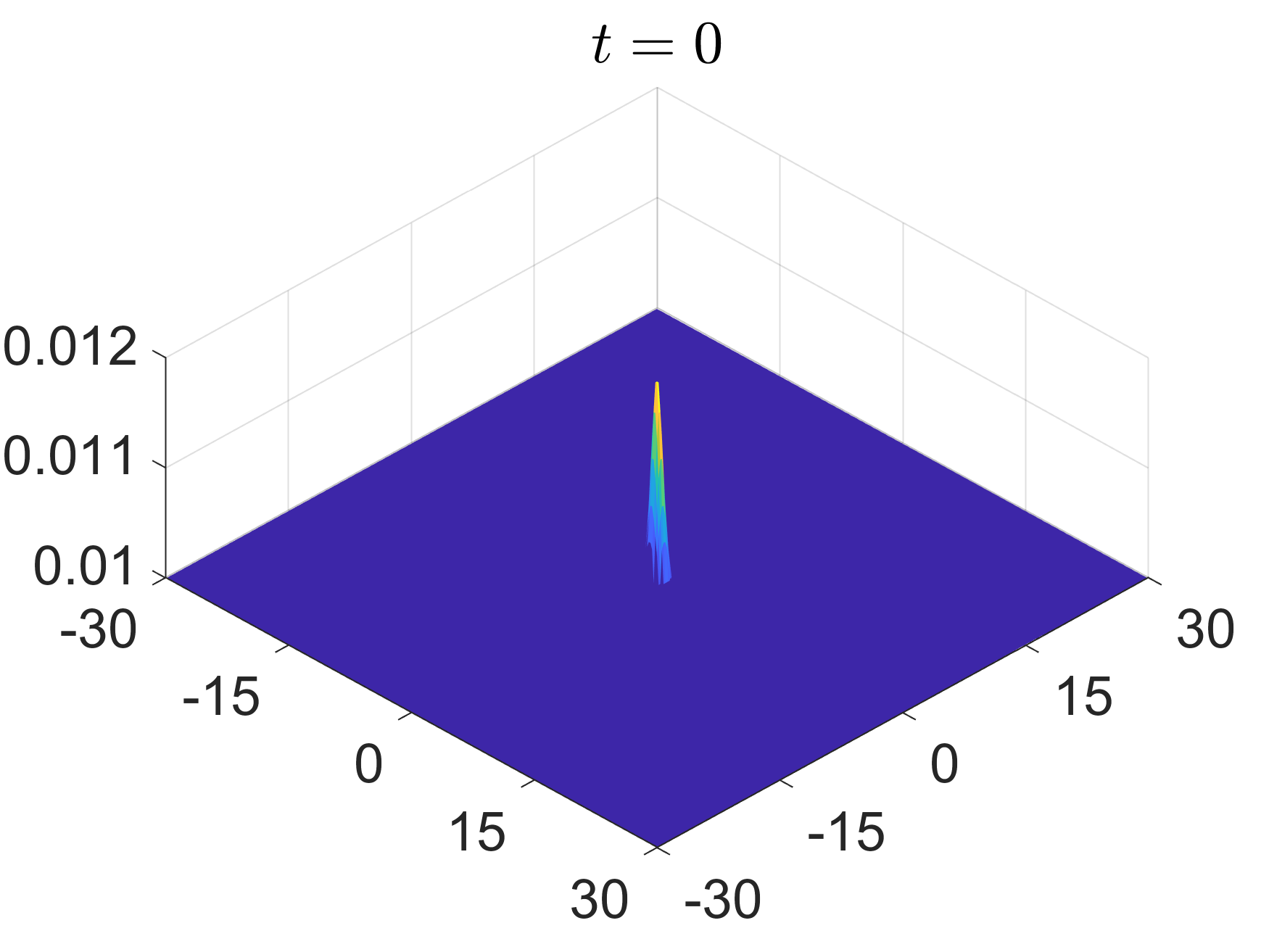}\hspace{1cm}
            \includegraphics[trim=0.0cm 0.3cm 0.3cm 0.2cm, clip, width=7.cm]{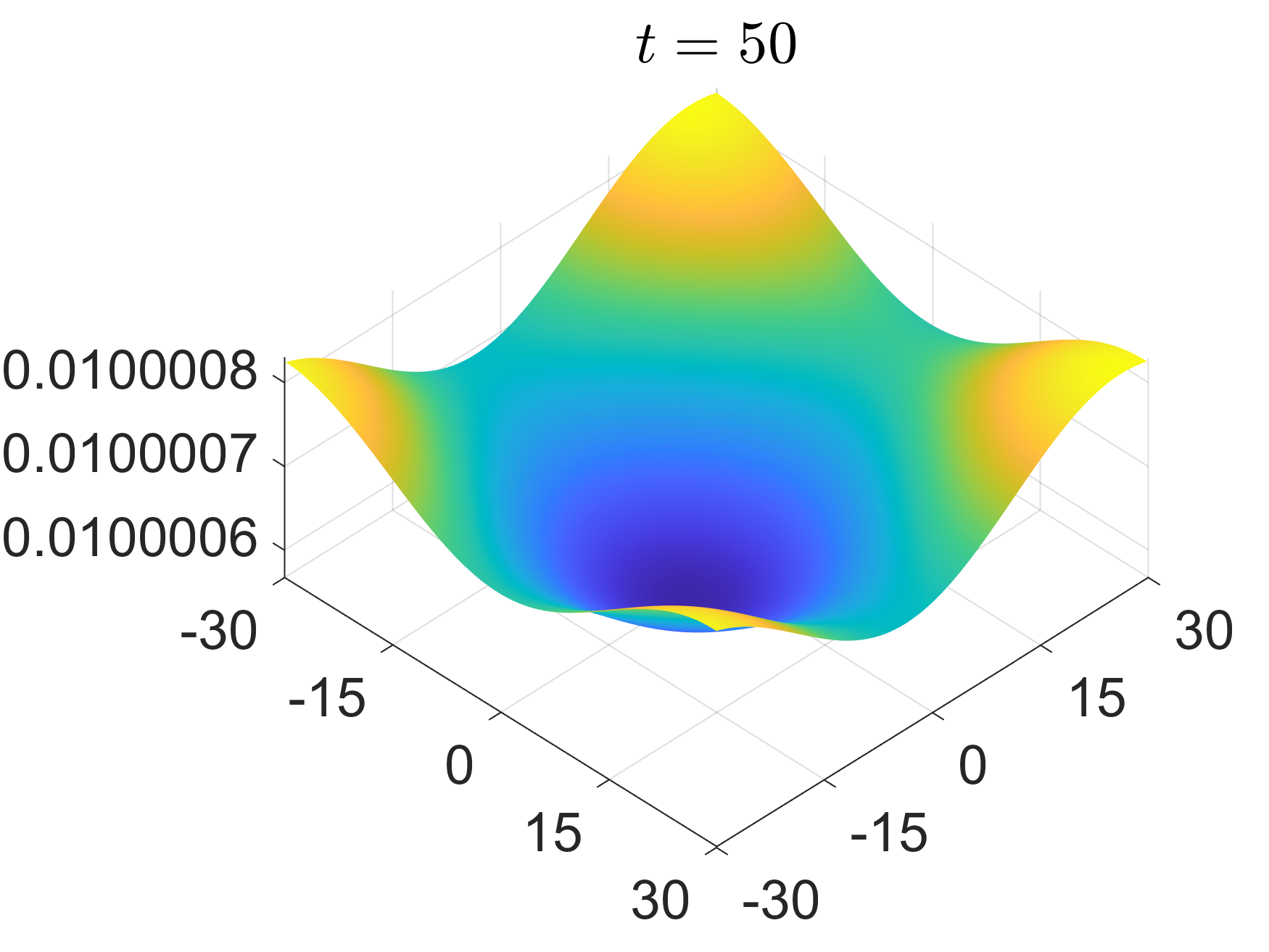}}
\vskip20pt
\centerline{\includegraphics[trim=0.0cm 0.3cm 0.3cm 0.2cm, clip, width=7.cm]{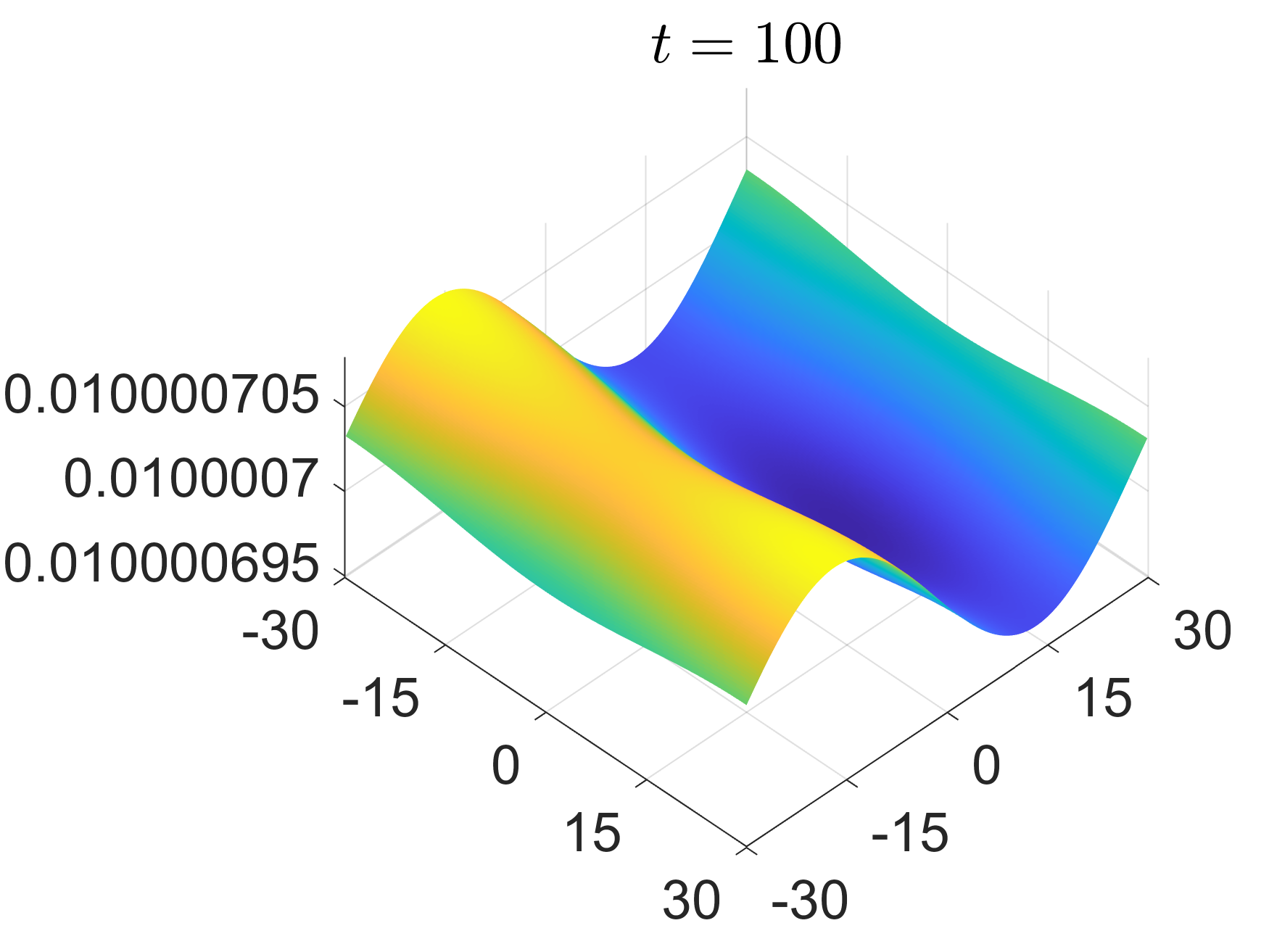}\hspace{1cm}
            \includegraphics[trim=0.0cm 0.3cm 0.3cm 0.2cm, clip, width=7.cm]{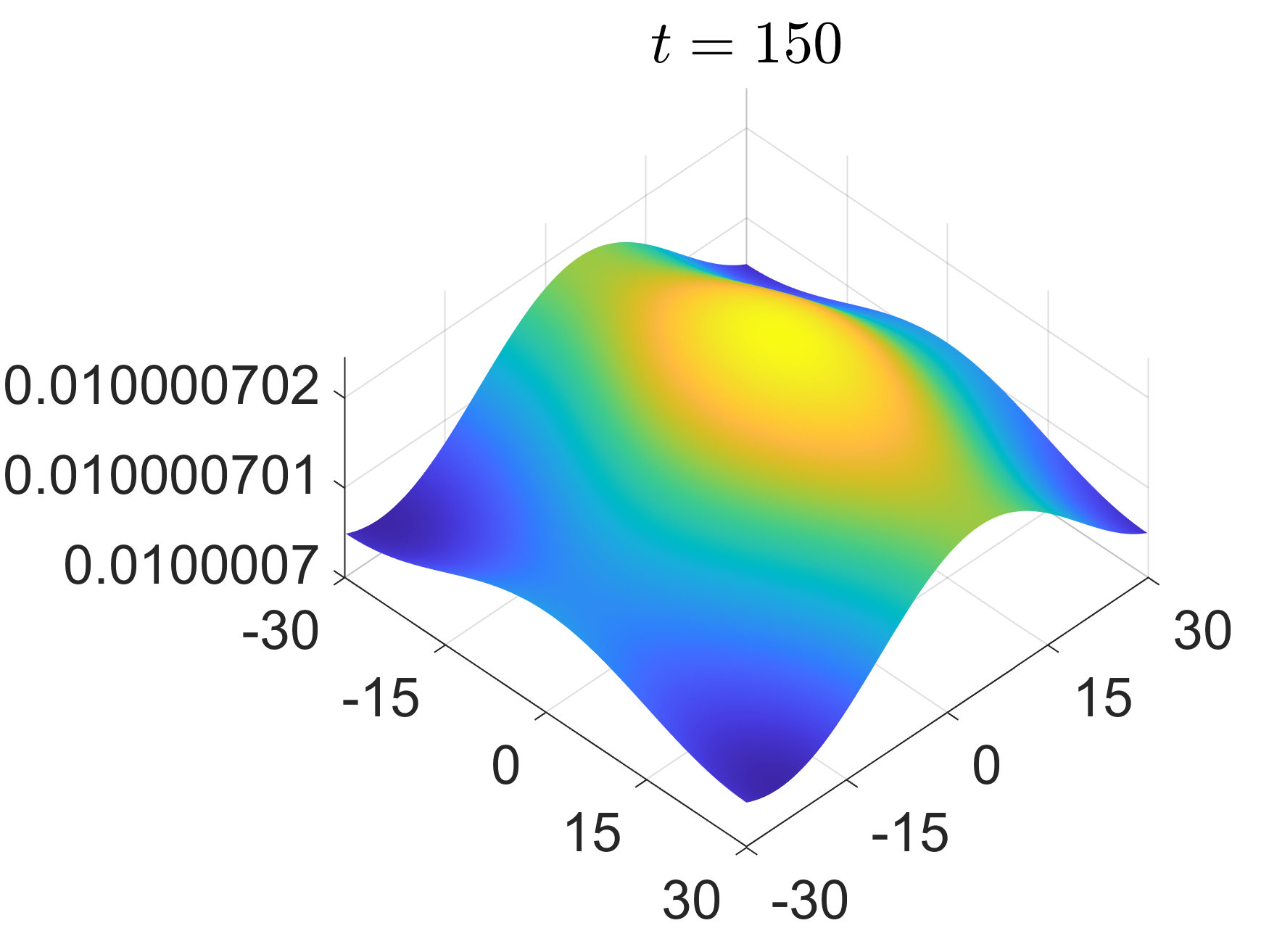}}
%\vskip20pt
%\centerline{\includegraphics[trim=0.0cm 0.3cm 0.3cm 0.2cm, clip, width=7.cm]{Example11_Fifth_h200}\hspace{1cm}
            %\includegraphics[trim=0.0cm 0.3cm 0.3cm 0.2cm, clip, width=7.cm]{Example11_Fifth_h250}}
\caption{\sf Linear case: $h$  computed with $\varepsilon=5$ (stable regime) at different times. \label{fig4}}
\end{figure}

\begin{figure}[ht!]
\centerline{\includegraphics[trim=0.0cm 0.3cm 0.3cm 0.2cm, clip, width=7.cm]{Initial_h0}\hspace{1cm}
            \includegraphics[trim=0.0cm 0.3cm 0.3cm 0.2cm, clip, width=7.cm]{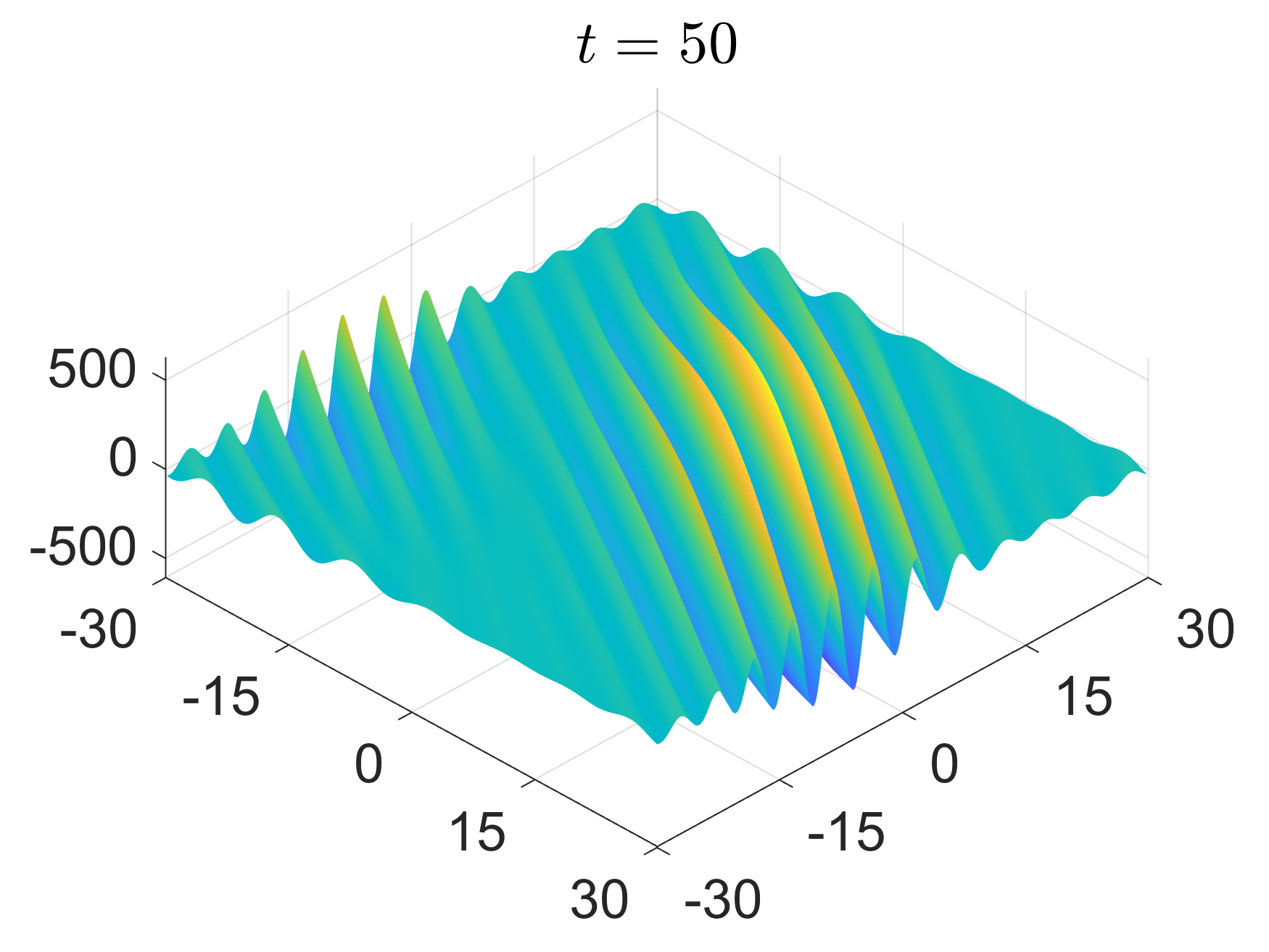}}
\vskip20pt
\centerline{\includegraphics[trim=0.0cm 0.3cm 0.3cm 0.2cm, clip, width=7.cm]{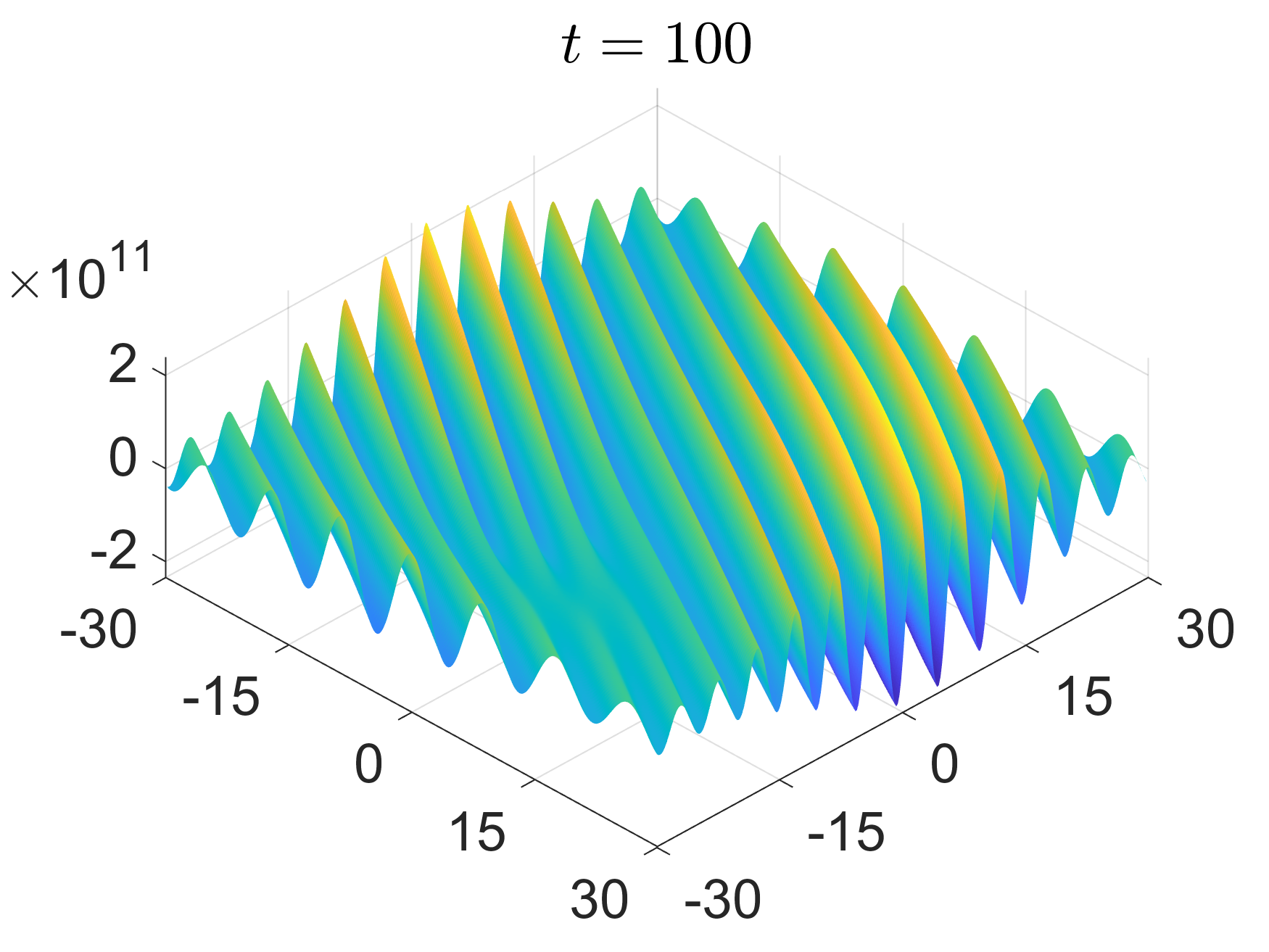}\hspace{1cm}
            \includegraphics[trim=0.0cm 0.3cm 0.3cm 0.2cm, clip, width=7.cm]{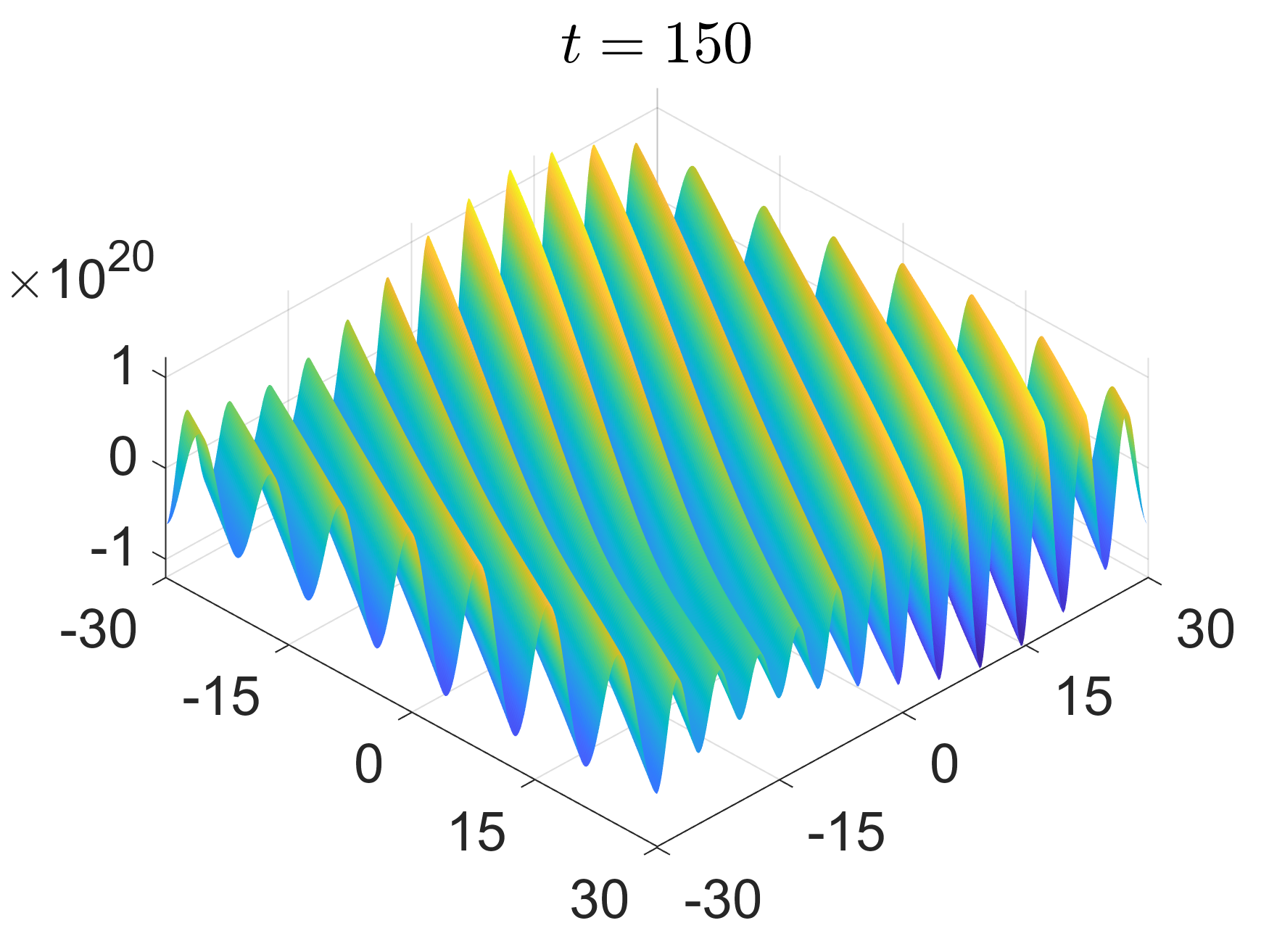}}
\vskip20pt
\centerline{\includegraphics[trim=0.0cm 0.3cm 0.3cm 0.2cm, clip, width=7.cm]{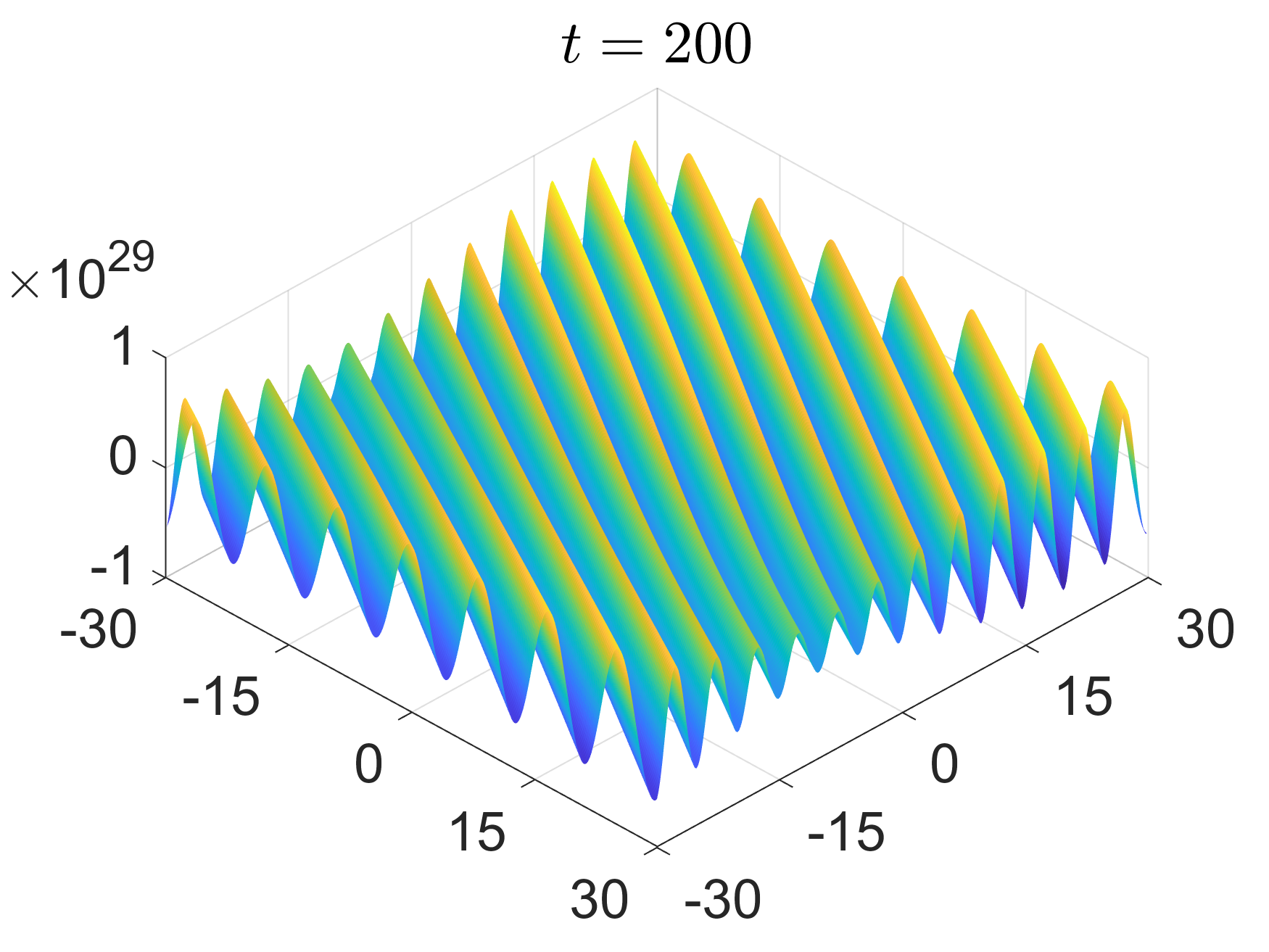}\hspace{1cm}
            \includegraphics[trim=0.0cm 0.3cm 0.3cm 0.2cm, clip, width=7.cm]{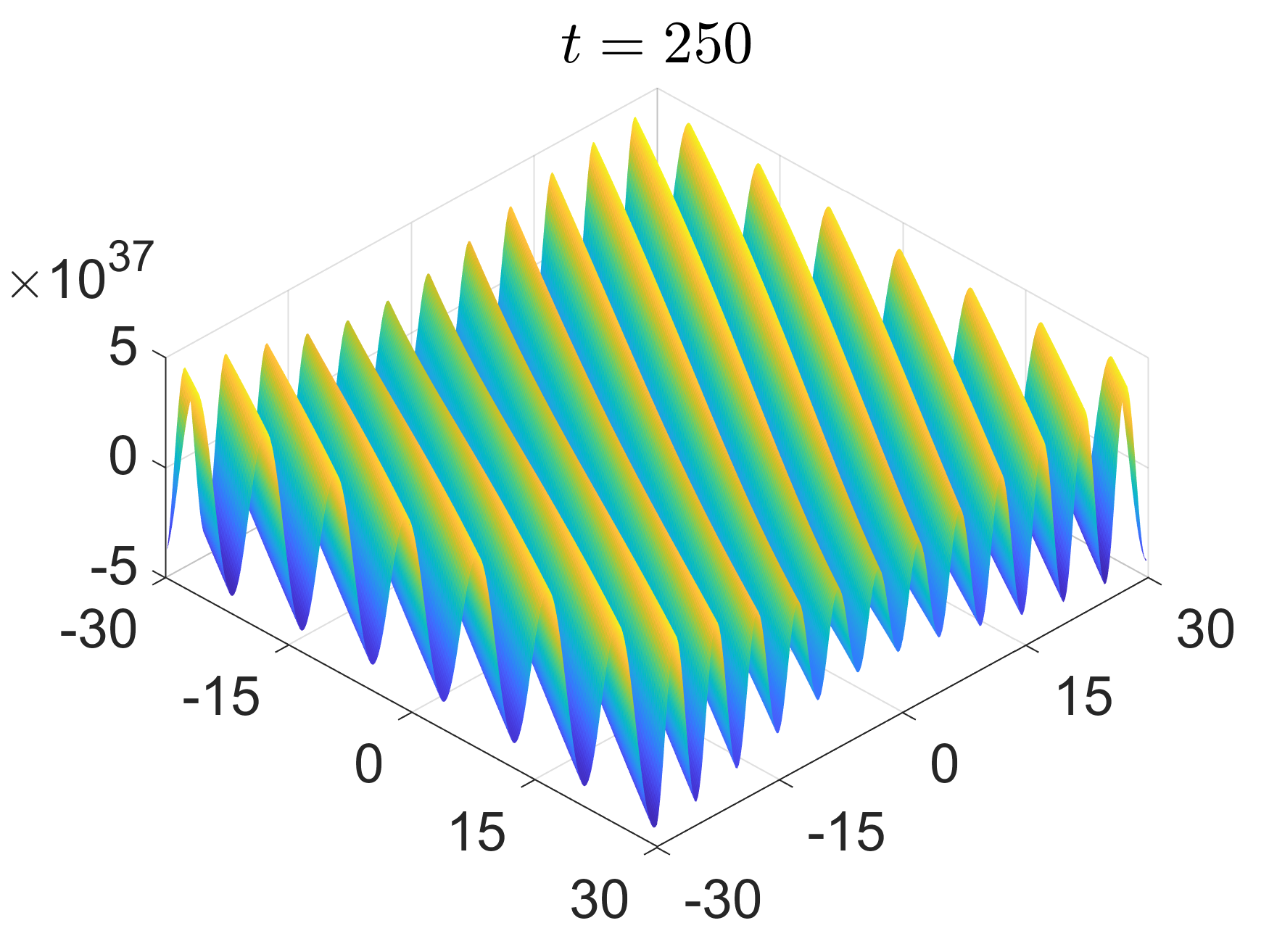}}
\caption{\sf Same as in Figure \ref{fig4}, but for $\varepsilon=1$ (unstable regime). \label{fig3}}
\end{figure}

In Figures \ref{fig4} and \ref{fig3}, we present the computed solutions for $\varepsilon=5$ and 1, respectively. One can clearly see that
the obtained numerical solutions support the analytical results in \S\ref{sec3}. Namely, when $\varepsilon=5$, the solution is stable and
the small initial perturbation decays in time as expected; see Figure \ref{fig4}, where we plot the computed solution at times $t=0$, 50, 
100, and 150. When $\varepsilon=1$, the solution is unstable: It creates large magnitude wave structures in the oblique direction that
corresponds to $\gamma=1/2$, and the magnitude of these waves grows indefinitely in time; see Figure \ref{fig3}, where we plot the computed
solution at times $t=0$, 50, 100, 150, 200, and 250.

\subsection{Nonlinear Case}
In this section, we conduct numerical experiments for the nonlinear viscous Saint-Venant system
\begin{equation*}
\bm u_t+\bm f(\bm U)_x+\bm g(\bm u)_y=C\bm u_{xx}+D\bm u_{yy}.
\end{equation*}
where
$$
\bm u=\begin{pmatrix}h\\q\\p\end{pmatrix},\quad
\bm f(\bm u)=\begin{pmatrix}q\\\dfrac{q^2}{h}+\dfrac{g}{2}h^2\\[1ex]\dfrac{pq}{h}\end{pmatrix},\quad
\bm g(\bm u)=\begin{pmatrix}p\\[0.5ex]\dfrac{pq}{h}\\[0.5ex]\dfrac{p^2}{h}+\dfrac{g}{2}h^2\end{pmatrix},
$$
and the matrices $C$ and $D$ are given by \eref{3.1}. We take the following initial conditions:
\begin{equation*}
h(x,y,0)=\begin{cases}1+\dfrac{1}{4000}e^{2(1-x^2-y^2)}&\mbox{if }x^2+y^2<1,\\1,&\mbox{otherwise},\end{cases}\quad
q(x,y,0)=p(x,y,0)\equiv0,
\end{equation*}
which are a small perturbation of the steady state $(h,q,p)\equiv(1,0,0)$.

\begin{figure}[ht!]
\centerline{\includegraphics[trim=0.0cm 0.3cm 0.3cm 0.2cm, clip, width=7.cm]{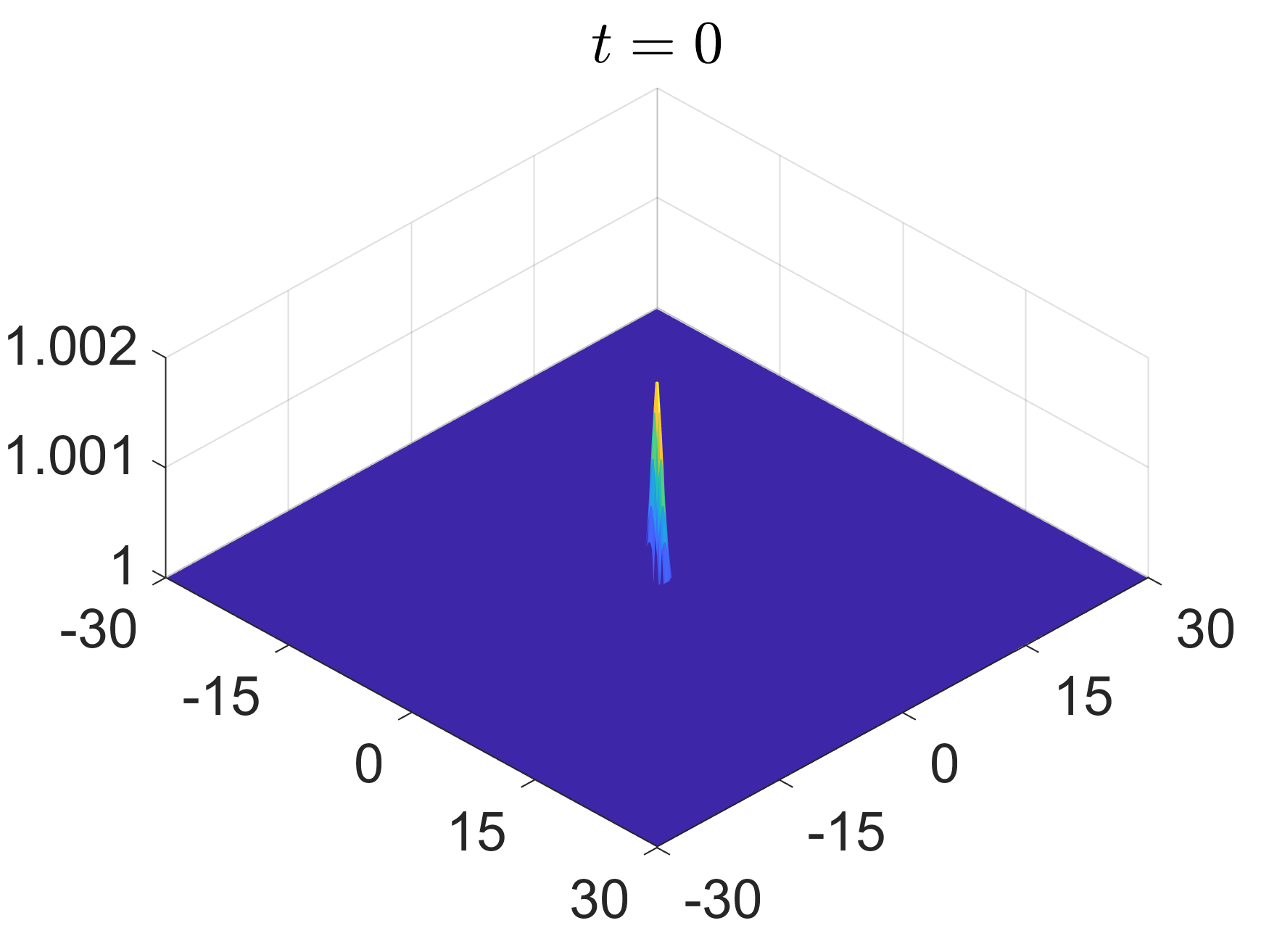}\hspace{0.5cm}
            \includegraphics[trim=0.0cm 0.3cm 0.3cm 0.2cm, clip, width=7.cm]{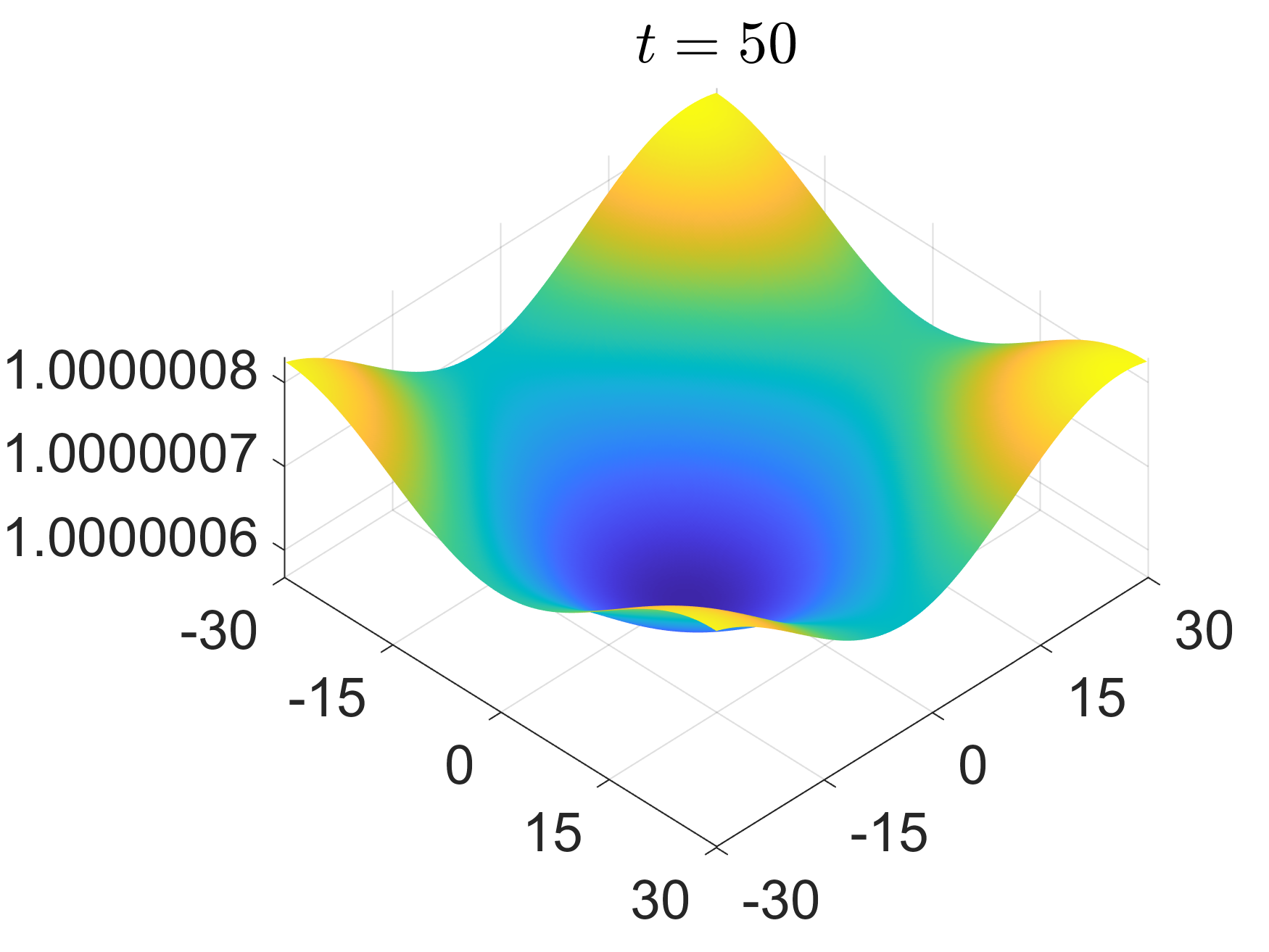}}
\vskip20pt
\centerline{\includegraphics[trim=0.0cm 0.3cm 0.3cm 0.2cm, clip, width=7.cm]{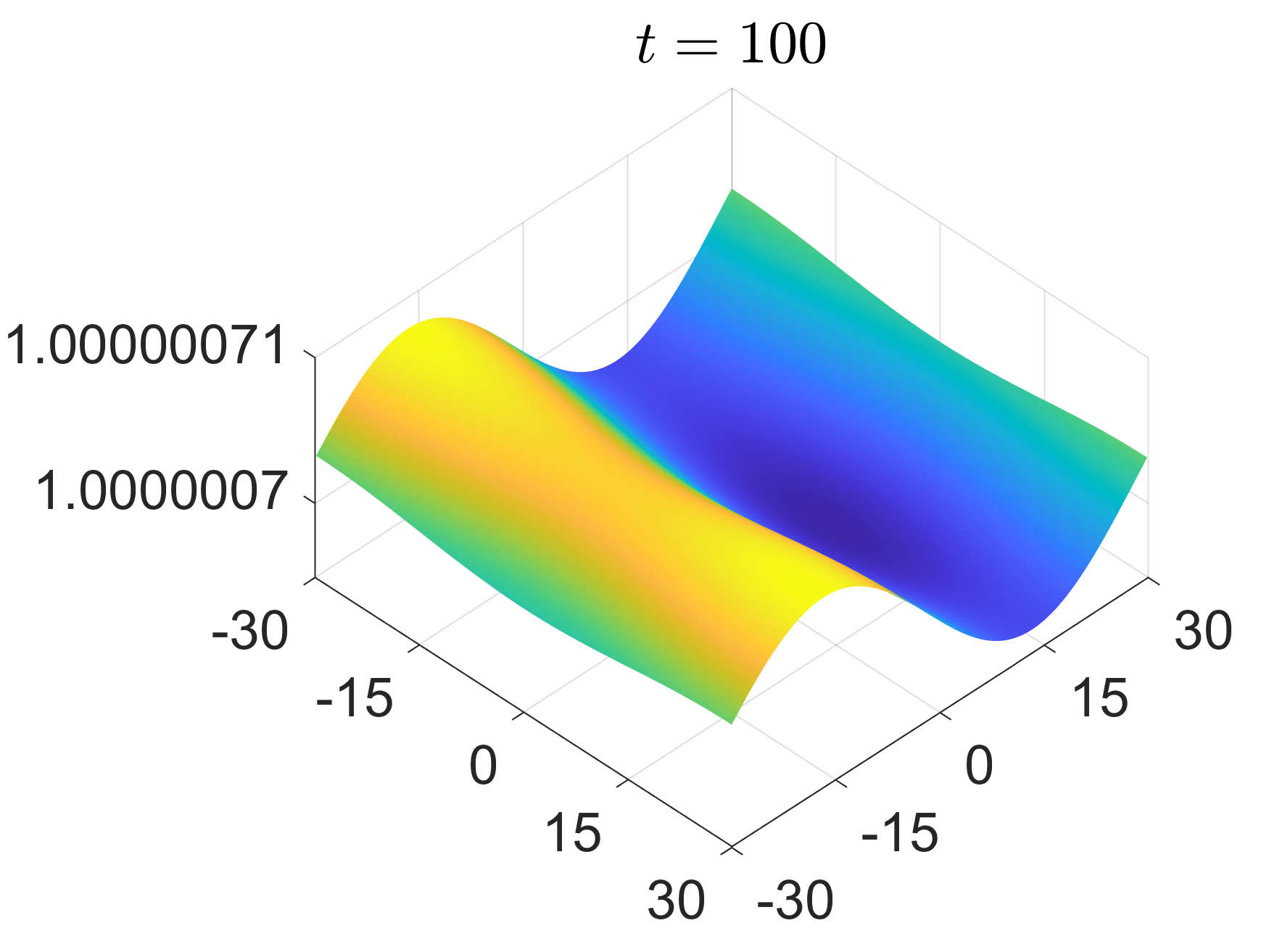}\hspace{1cm}
            \includegraphics[trim=0.0cm 0.3cm 0.3cm 0.2cm, clip, width=7.cm]{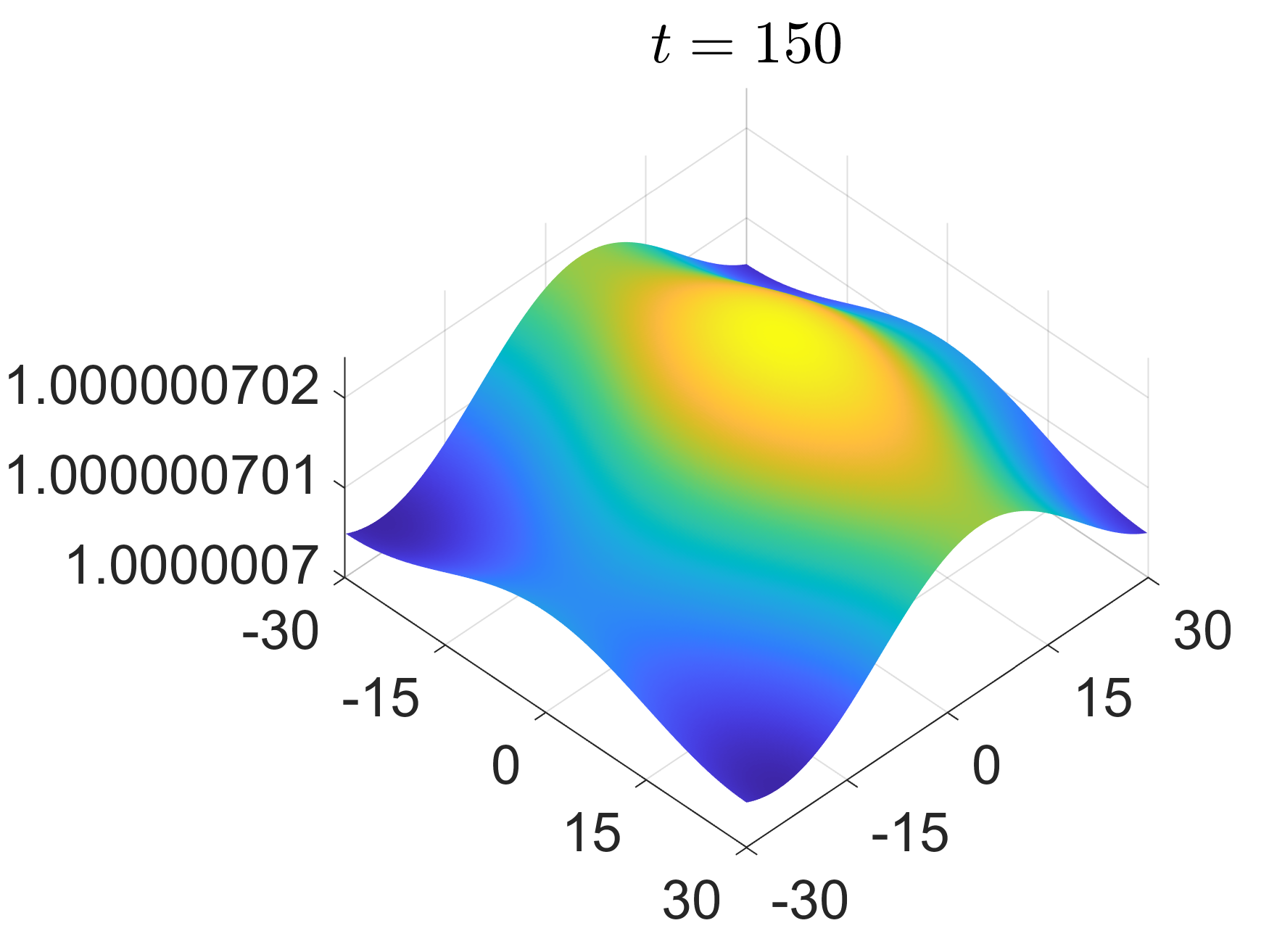}}
%\vskip20pt
%\centerline{\includegraphics[trim=0.0cm 0.3cm 0.3cm 0.2cm, clip, width=7.cm]{Example2_Fifth_h200}\hspace{1cm}
            %\includegraphics[trim=0.0cm 0.3cm 0.3cm 0.2cm, clip, width=7.cm]{Example2_Fifth_h250}}
\caption{\sf  Nonlinear case: $h$  computed with $\varepsilon=5$ (stable regime) at different times. \label{fig6}}
\end{figure}

\begin{figure}[ht!]
\centerline{\includegraphics[trim=0.0cm 0.3cm 0.3cm 0.2cm, clip, width=7.cm]{Initial_h}\hspace{0.5cm}
            \includegraphics[trim=0.0cm 0.3cm 0.3cm 0.2cm, clip, width=7.cm]{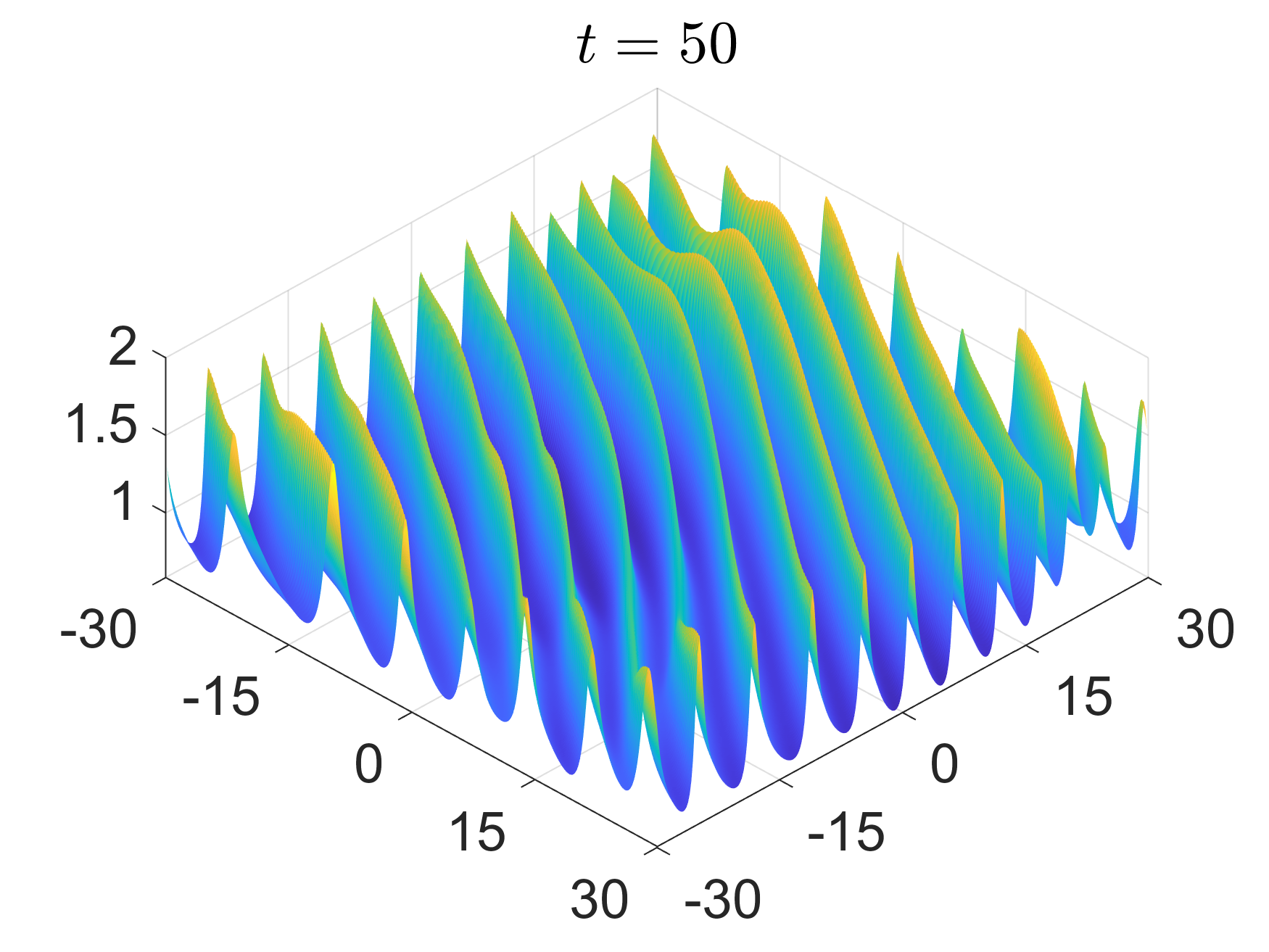}}
\vskip20pt
\centerline{\includegraphics[trim=0.0cm 0.3cm 0.3cm 0.2cm, clip, width=7.cm]{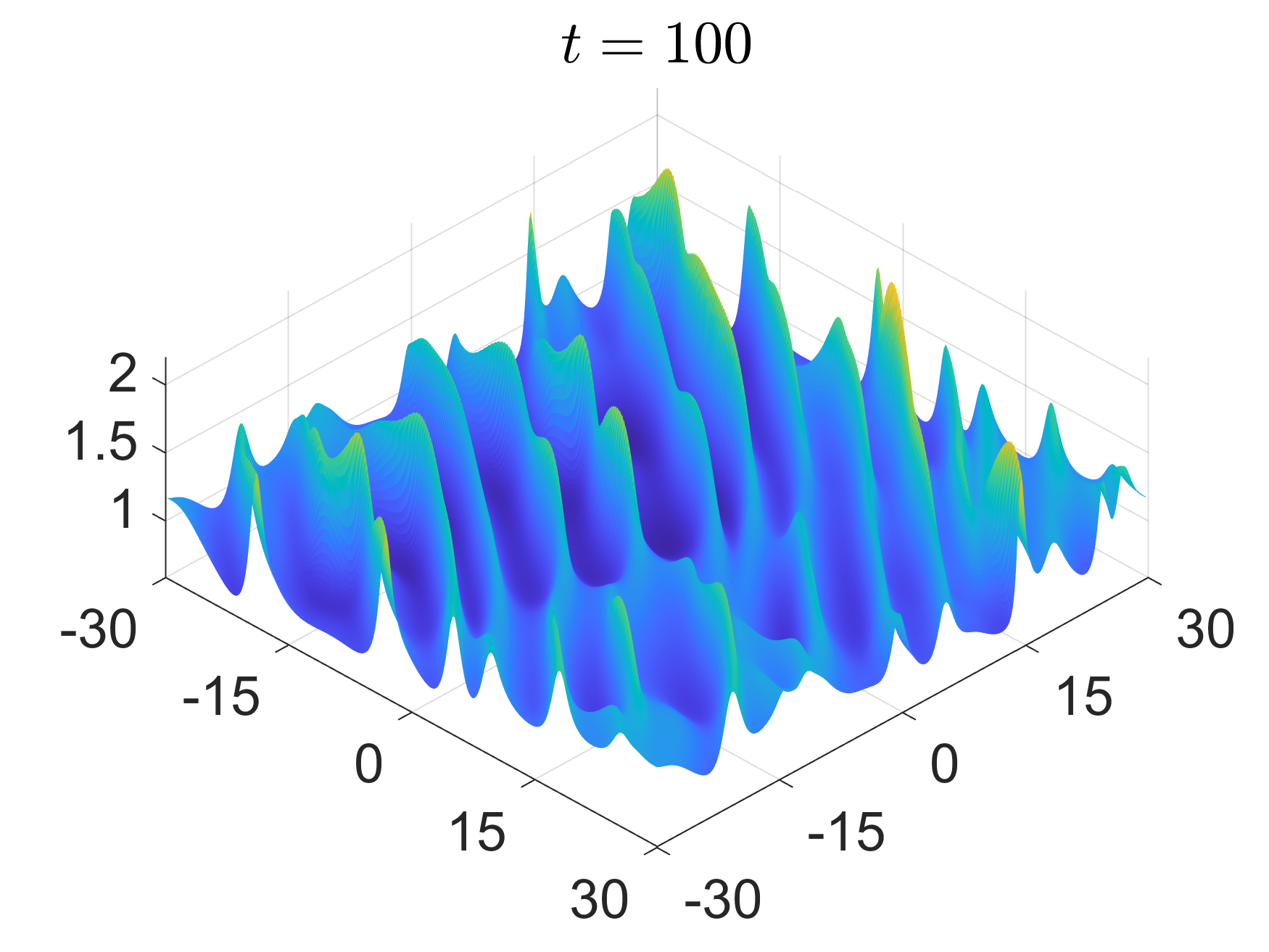}\hspace{1cm}
            \includegraphics[trim=0.0cm 0.3cm 0.3cm 0.2cm, clip, width=7.cm]{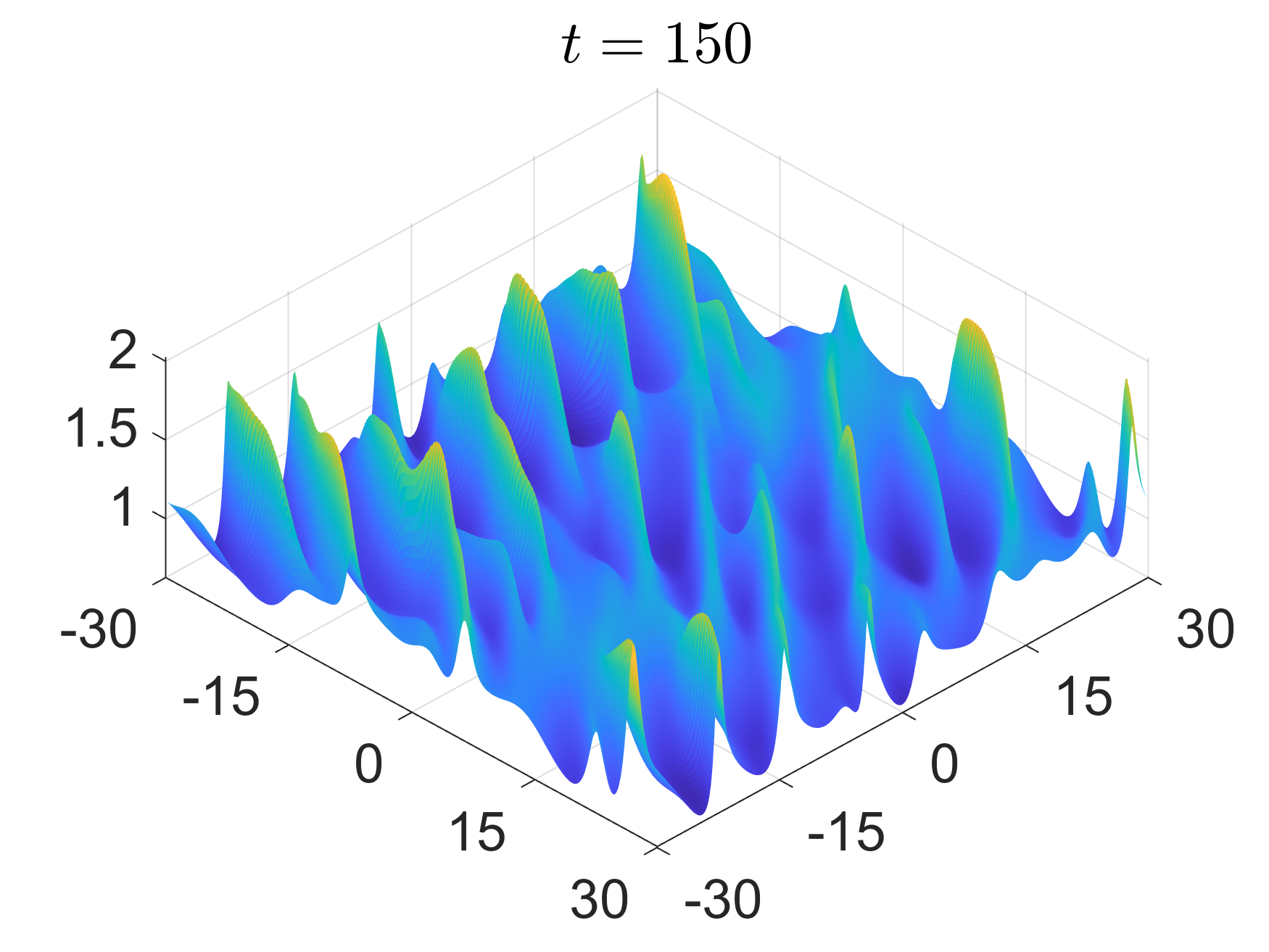}}
\vskip20pt
\centerline{\includegraphics[trim=0.0cm 0.3cm 0.3cm 0.2cm, clip, width=7.cm]{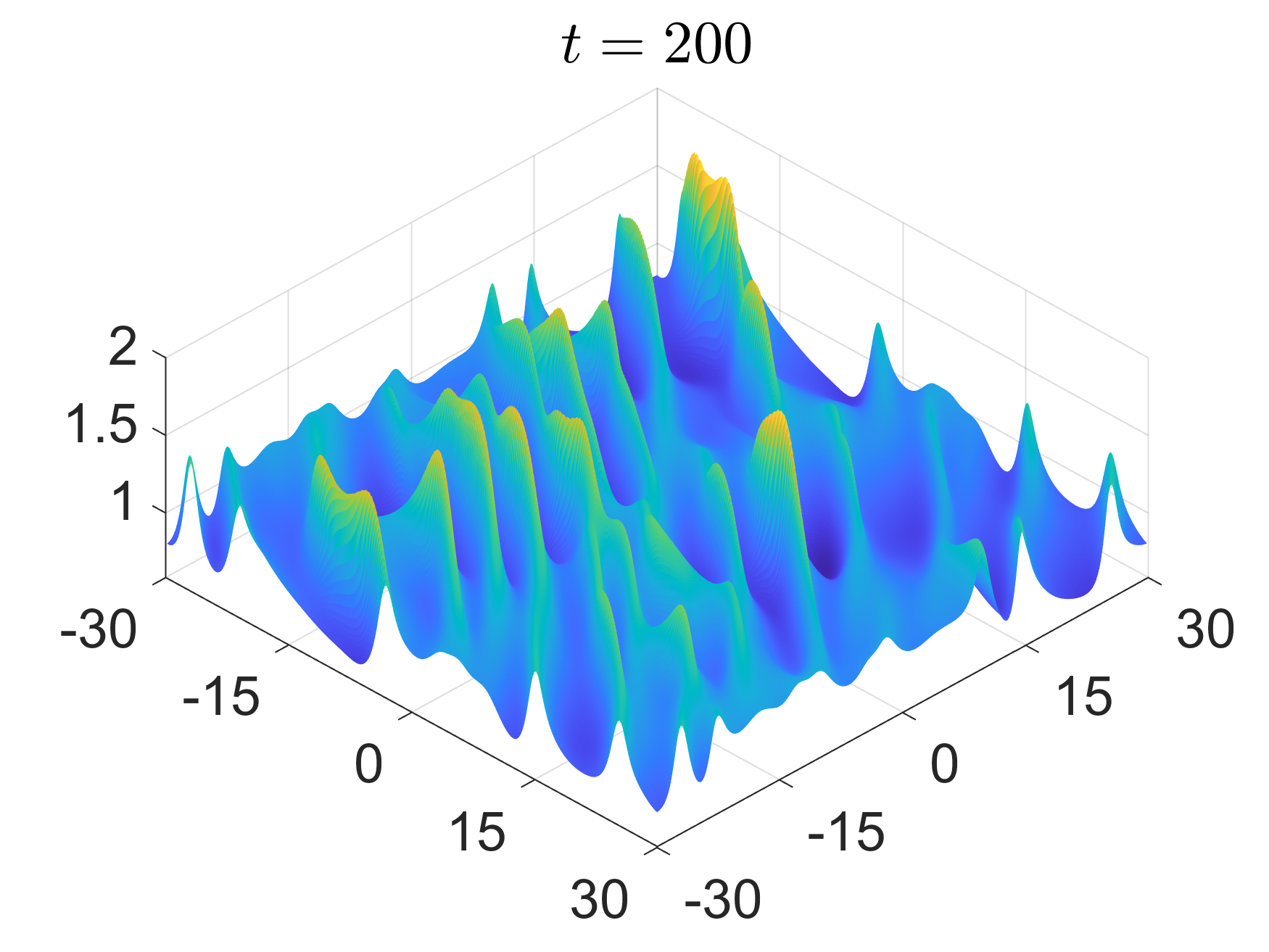}\hspace{1cm}
            \includegraphics[trim=0.0cm 0.3cm 0.3cm 0.2cm, clip, width=7.cm]{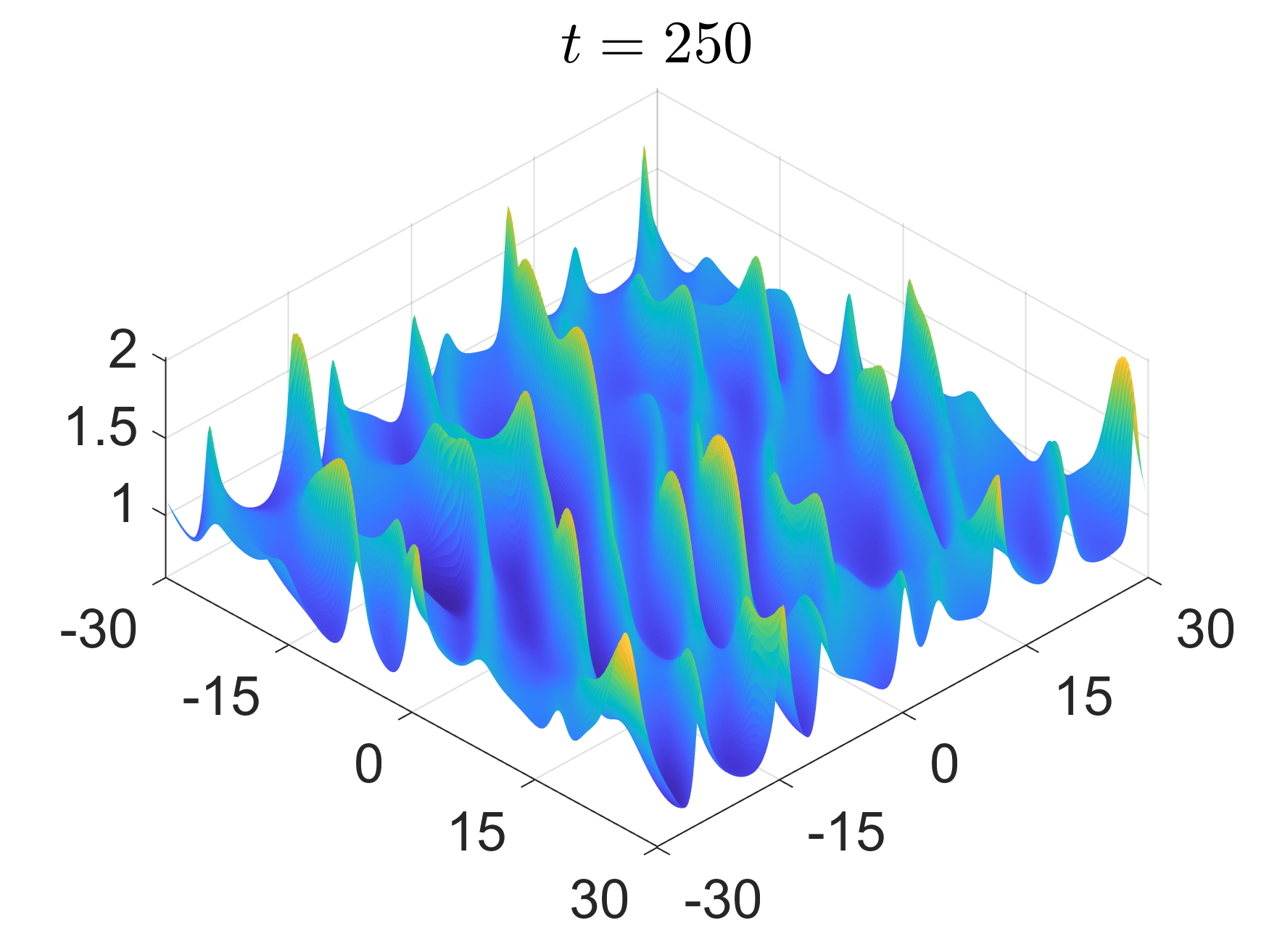}}
\caption{\sf Same as in Figure \ref{fig6}, but for $\varepsilon=1$ (unstable regime). \label{fig5}}
\end{figure}

We present the computed solutions in Figures \ref{fig6} and \ref{fig5} for $\varepsilon=5$ and 1, respectively. As one can see, when
$\varepsilon=5$ the solution of the nonlinear system is stable and, like in the linear case, the initial perturbation decays in time; see
Figure \ref{fig6}, where we plot the computed solution at times $t=0$, 50, 100, and 150. The solution behavior in the unstable
($\varepsilon=1$) regime is, however, different from the linear case. While the solution develops large magnitude waves in the same oblique
direction corresponding to $\gamma=1/2$, the magnitude of these waves does not increase beyond 2 and the solution structure evolves in a
rather complicated nonlinear manner; see Figure \ref{fig5}, where we plot the computed solution at times $t=0$, 50, 100, 150, 200, and 250.

\section{Stability Conjectures}
In this section, we discuss a general question on conditions on the matrices $A$, $B$, $C$, and $D$ for the stability of the system
\eref{1.2}. Those ideas are prompted by simple ``combinatorial'' observations of 1-D systems in the $x$- and $y$-directions, as well as by
the examples of unstable oblique waves above. We surmize that linear operators, stemming from all combinations of convective terms in all
directions with viscous terms in all directions, should be stable for all wavenumbers. 
\begin{conj}\label{con51}
If the eigenvalues of the stability operators $\Omega_1(k)=-ikA-k^2C$, $\Omega_2(k)=-ikA-k^2D$, $\Omega_3(k)=-ikB-k^2C$, and
$\Omega_4(k)=-ikB-k^2D$ are negative for all $k$, then the system \eref{1.2} is stable.
\end{conj}
\begin{rmk}
In fact, Conjecture \ref{con51} suggests to check the eigenvalues of four possible combinations $-ik(\cdot)-k^2(\cdot)$ of all convective
and viscous matrices. 
\end{rmk}

Before presenting the second conjecture, we will make a few simple observations. We first consider a general eigenvalue problem prompted by
the stability considerations. The eigenvalues of $A$ are the roots of the characteristic polynomial of degree $N$,
$P_N(\lambda)=\sum_{i=0}^N(-1)^ir_i\lambda^{N-i}$ with $r_0=1$, $r_1={\rm tr}(A)$, and $r_N=\det(A)$. In fact, the coefficients $r_i$ for
$i=1,\ldots,N$ are the sum of the $i$-rowed principal minors of $A$. In particular, for $N=2$ and $N=3$, we have
\begin{equation*}
\begin{aligned}
&P_2(\lambda)=\lambda^2-{\rm tr}(A)\lambda+\det(A)=\lambda^2-(a_{11}+a_{22})\lambda+(a_{11}a_{22}-a_{12}a_{21}),\quad
P'_2(\lambda)=2\lambda-{\rm tr}(A),\\
&P_3(\lambda)=\lambda^3-{\rm tr}(A)\lambda^2+\Big(\sum_{i=1}^3M_{ii}\Big)\lambda-\det(A),\quad
P'_3(\lambda)=3\lambda^2-2\,{\rm tr}(A)\lambda+\sum_{i=1}^3M_{ii}.
\end{aligned}
\end{equation*}
Here, $M_{ij}$ denotes the $(i,j)$ minor of $A$, and
$\sum_{i=1}^3M_{ii}=(a_{11}a_{22}-a_{12}a_{21})+(a_{11}a_{33}-a_{13}a_{31})+(a_{22}a_{33}-a_{23}a_{32})$. 

Let us denote by $\lambda_\ell$, $\ell=1,\ldots,N$ the eigenvalues of $A$ and consider a slightly {\em perturbed} matrix $A+\delta C$, where
$\delta$ is a small parameter. Its eigenvalues are $\lambda_\ell+\delta\lambda_{\ell,1}+\delta^2\lambda_{\ell,2}+\ldots$ and the
characteristic equation can be written as
\begin{equation}
\widetilde P_N(\lambda_\ell;\delta):=\sum_{i=0}^N(-1)^i(r_i+\delta r_{i,1}+\delta^2r_{i,2}+\ldots)
(\lambda_\ell+\delta\lambda_{\ell,1}+\delta^2\lambda_{\ell,2}+\ldots)^{N-i}=0
\label{5.1}
\end{equation}
with the coefficients on the left-hand side being slightly perturbed sums of all principal minors of $A+\delta C$. Differentiating
\eref{5.1} with respect to $\delta$ and substituting $\delta=0$ results in
\begin{equation*}
\frac{\partial\widetilde P_N}{\partial\delta}(\delta=0)=\lambda_{\ell,1}P'_N(\lambda_{\ell})+\sum_{i=1}^N(-1)^ir_{i,1}\lambda_\ell^{N-i}=0,
\end{equation*}
which implies that the linear correction of the eigenvalue $\lambda_\ell$, which may influence the sign of the $\ell$-th eigenvalue of the
perturbed matrix for small $\delta$, is
\begin{equation}
\lambda_{\ell,1}=-\frac{1}{P'_N(\lambda_{\ell})}\sum_{i=1}^N(-1)^ir_{i,1}\lambda_\ell^{N-i},\quad\ell=1,\ldots,N.
\label{5.2}
\end{equation}
Here, $r_{i,1}$ are the sum of the $i$-rowed principal minors with all combinations of elements where precisely one element of matrix $A$ is
replaced by elements of matrix $C$ (first-order correction) so that $\lambda_{\ell,1}$ is a linear function relative to the elements of
matrix $C$. Notice that $P'_N(\lambda_{\ell})\neq0$ as long as the eigenvalues of $A$ are distinct.

In particular, for $N=2$, the first correction \eref{5.2} assumes a simple elegant form:
\begin{equation*}
\lambda_{\ell,1}=\frac{\left(c_{11}+c_{22}\right)\lambda_\ell+a_{12}c_{21}-a_{22}c_{11}+a_{21}c_{12}-a_{11}c_{22}}
{2\lambda_\ell-(a_{11}+a_{22})},\quad\ell=1,2. 
\end{equation*}
For $N=3$, \eref{5.2} takes the following (a bit cumbersome) form:
\begin{equation*}
\lambda_{\ell,1}=\frac{{\rm tr}(C)\lambda_\ell^2-r_{2,1}\lambda_\ell+\displaystyle{\sum_{i=1}^3\sum_{j=1}^3(-1)^{i+j}c_{ij}M_{ij}}}
{3\lambda^2_\ell-2\,{\rm tr}(A)\lambda_\ell+\displaystyle{\sum_{i=1}^3M_{ii}}},\quad\ell=1,2,3,
\end{equation*}
where
$$
r_{2,1}=a_{22}c_{33}+a_{33}c_{22}-a_{23}c_{32}-a_{32}c_{23}+a_{11}c_{33}+a_{33}c_{11}-a_{13}c_{31}-a_{31}c_{13}+
a_{11}c_{22}+a_{22}c_{11}-a_{12}c_{21}-a_{21}c_{12}
$$
 
We now use \eref{5.2} and introduce the second conjecture, which relates to a more algebraic reduction of Conjecture \ref{con51}.
\begin{conj}\label{con53}
Let the matrix $A$ has distinct real eigenvalues $\lambda_\ell,~\ell=1,\ldots,N$, and the matrix $C$ has real positive eigenvalues. Then,
the linear operator $\Omega(k)=-ikA-k^2C$ will have negative eigenvalues if and only if the quantities $\lambda_{\ell,1}$ defined in
\eref{5.2} are negative for all $\ell=1,\ldots,N$.
\end{conj}

We stress that values $\lambda_{\ell,1}$ control the first correction of long-wavelength expansion for small $k$:
$-i\lambda_\ell k+\lambda_{\ell,1}k^2$. The importance of Conjecture \ref{con53} is that it gives a simple algebraic stability criterion
based on the elements of matrices $A$ and $C$ for {\em all} $k$. Both Conjectures \ref{con51} and \ref{con53} may be viewed as a direct and
simple solution of the multidimensional Gelfand problem, and they are inspired by the remarkable paper by Majda and Pego \cite{MP85}.
 
\section{Conclusion}
We have demonstrated the counter-intuitive extension of Majda-Pego instability for 2-D systems, both analytically and through direct
numerical simulations. To the best of our knowledge, this is the first nontrivial result for the Gelfand problem in the
multidimensional case.

While we have presented examples of unstable 2-D system with viscous terms, it is rather easy to construct similar examples in the 3-D case:
as in the 2-D case, both advective and viscous matrices should be non-diagonal, all of the waves propagating in either $x$-, $y$-, or
$z$-direction will be stable, yet some oblique waves will be unstable.

We hope that our examples will stimulate the elucidation of peculiar role of viscous effects in the broad spectrum of natural phenomena,
including fluid mechanics. We have proposed the conjectures that give direct assessment of stability of multidimensional
convection-diffusion systems and resolve the Gelfand problem. One can also generalize Conjectures \ref{con51} and \ref{con53} for the 3-D
systems: In order to establish the stability criteria, one has to study to the eigenvalues of nine linear operators $-ik(\cdot)-k^2(\cdot)$
for possible combination of all convective and viscous matrices.

\begin{acknowledgment}
The work of A. Kurganov was supported in part by NSFC grant 12171226 and the fund of the Guangdong Provincial Key Laboratory of
Computational Science and Material Design (No. 2019B030301001).
\end{acknowledgment}

\bibliography{ref}
\bibliographystyle{siam}
\end{document}